\documentstyle[11pt]{article}
\newtheorem {th}{Theorem}[section]
\newtheorem {pr}[th]{Proposition}
\newtheorem {lem}[th]{Lemma}
\newtheorem {cor}[th]{Corollary}

\def\Cox{\hfill \Box}

\def\sf{\sigma\mbox{-field}}
\def\dd{\delta}
\def\ee{\epsilon}
\def\E{{\bf{E}}}
\def\P{{\bf{P}}}
\def\N{\hbox{I\kern-.2em\hbox{N}}}
\def\R{\hbox{I\kern-.2em\hbox{R}}}

\def\Z{{\bf{Z}}}
\def\F{{\cal{F}}}

\def\|{\, | \, }
\def\v0{{\bf 0}}
\def\one{{\bf 1}}
\def\0{\hat{0}}
\def\1{\hat{1}}

\def\Bt{{\tilde{B}}}

\begin{document}
\begin{center}
{\large \bf RANDOM WALKS IN VARYING DIMENSIONS} \\
\end{center}
\vspace{2ex}
\begin{center}
{\sc Itai Benjamini} \footnote{Mathematical Sciences Institute, 409 College
Ave., Ithaca, NY, 14853. Research partially supported by the U. S.
Army Research Office through the Mathematical Sciences Institute
of Cornell University.}, \, 
{\sc Robin Pemantle} \footnote{Research supported in part by 
National Science Foundation grant \# DMS 9300191, by a Sloan Foundation
Fellowship, and by a Presidential Faculty Fellowship.}$^,$\footnote{Department 
of Mathematics, University of Wisconsin-Madison, Van Vleck Hall, 480 Lincoln
Drive, Madison, WI 53706 \, .}, \,  and \, 
{\sc Yuval Peres} \footnote{Department of
Statistics, 367 Evans Hall University of California, Berkeley, CA 94720.
Research partially supported by NSF grant
\# DMS-9404391 and by a Junior Faculty Fellowship from the Regents of
the University of California.}
\end{center}
\vspace{5ex}
\begin{center}
{\bf Abstract}
\end{center}
We establish  recurrence criteria for sums of independent random 
variables which take values in Euclidean lattices of varying dimension.
 In particular, we describe  transient inhomogenous random walks
in the plane which interlace two symmetric step distributions of bounded
support.

\newpage

\section{Introduction}

As everyone knows, and Polya proved in 1921,
 simple random walk is recurrent in one
and two dimensions and transient in three or more dimensions.
Also widely known is that the transition between recurrence
and transience occurs precisely at dimension 2,
rather than at some fractional dimension in the interval $(2,3)$.
  This is not
a precise statement, but one common interpretation is:
\begin{quote}
When the dimension parameter $d$ in formulae for quantities 
such as  Green's function is taken to be continuous, 
then qualitatively different behavior is observed in
the two regimes $d \leq 2$ and $d > 2$.  
\end{quote}
One can interpret the phase transition 
at dimension $2$ probabilistically by
considering simple random walk on subgraphs of the three-dimensional
lattice; this was done by T. Lyons (1983), who showed
that a ``slight fattening of a quadrant in $\Z^2$''
suffices to obtain transience.  Another approach, taken here,
 is to 
construct a random walk on a $d$-dimensional lattice which
only occasionally, at some fixed times, moves in directions outside a certain 
subspace of smaller dimension.  The recurrence/transience
of such a walk depends (as one would expect) on the frequency
of the fully $d$-dimensional steps, but the location of the
phase boundary is not what one would predict from a look
at the Green function; in other words, the usual ``Borel-Cantelli''
criterion for transience fails. A related model, in which the
exceptional moves are taken at random times, was analyzed by Scott (1990).

Let $F_3$ be a truly three
dimensional distribution with mean zero and finite variance on the lattice
$\Z^3$, and denote by $F_2$ the projection of this distribution
to the $x$-$y$ plane.  Given an increasing
sequence of positive integers $\{ a_n \}$, we consider the inhomogeneous
random walk $\{S_k\}$ whose independent increments $S_k - S_{k-1}$ 
have distribution $F_3$ if $k \in \{ a_n \}$ and distribution $F_2$
otherwise.  Theorem~\ref{th4.4} shows that the process $\{S_k\}$
is {\em recurrent} if 
$$a_n \approx \exp ( \exp ( n^{1/2} ))$$
but {\em transient} if $a_n \approx \exp ( \exp ( n^{\theta} ))$
for any $\theta \in (0,1/2)$.  Here {\em recurrence} means that the
number of $k$ for which $S_k = \v0$ is almost surely infinite,
and {\em transience} means that this number is almost surely finite.
(These alternatives are exhaustive, cf. Lemma~\ref{lem4.1}.)
An easy calculation shows that the {\em expected} number of visits
to the origin by $\{S_k\}$ is infinite when $0< \theta < 1/2$ as well 
as when $\theta = 1/2$.  

We also consider variants in other dimensions.  For instance,
there exists a recurrent random walk which interlaces two-dimensional,
four-dimensional and six-dimensional steps (but the four-dimensional
steps are indispensible here; see Corollary~\ref{cor4.3}).  
Conversely, there is a transient process obtained by alternating
blocks of one-dimensional and two-dimensional  random walk steps,
where both increment distributions are symmetric and have bounded support
(Proposition~\ref{pr6.1}).  

\section{Statement of results}

To give meaning to the terms ``recurrent'' and ``transient'', we 
first prove a ``folklore'' lemma 
which implies a 0-1 law for recurrence of RWVD.  

\begin{lem} \label{lem4.1}
Let $\{F_j : 1 \leq j \leq l \}$ be distributions on the abelian group
$Y$ and let 
$$(n(1) , n(2) , \ldots ) \in \{1 , 2 , \ldots , l \}^{\Z^+}$$
be any sequence in which each value $1 , \ldots , l$ occurs
infinitely often.  Let $\{ X_k \}$ be independent random variables
with corresponding distributions $\{F_{n(k)}\}$. 
Then any tail event for the sequence
of partial sums $S_N = \sum_{k=1}^N X_k$ has probability 0 or 1.
\end{lem}

\noindent{\sc Proof:}  If $l = 1$, this is a consequence of the
Hewitt-Savage 0-1 law.  If $l > 1$, assume for induction that 
the result is true for smaller values of $l$, and let $\F_{l-1}$
denote the $\sf$ generated by 
\begin{equation} \label{eq cond ex}
\{ X_k : n(k) \leq l-1 \} .  
\end{equation}
Conditional on $\F_{l-1}$, the event $B$ is exchangeable in
the remaining variables \newline $\{X_k : n(k) = l \}$; since these
variables are identically distributed, the Hewitt-Savage 0-1 law
shows that $\P [B \| \F_{l-1}] \in \{0 , 1 \}$ almost surely.
The set $\Bt := \{ \P [B \| \F_{l-1}] = 1 \}$ is 
$\F_{l-1}$-measurable, and it is a tail event for the partial sums of
the variables in~(\ref{eq cond ex}).  By induction, $\P [\Bt] \in \{0,1\}$, 
which shows that $\P [B] \in \{0,1\}$.    $\Cox$

\noindent{\bf Definition:} Let $d < D$ be positive integers and let
$\{a_n\}$ be an increasing sequence of integers.
A $\Z^d${\em in} $\Z^D$ {\em random walk in varying dimension } is a process
 $\{ S_k \}$ in $\Z^D$ with independent increments
$S_k - S_{k-1}$ distributed according to a truly $D$-dimensional
distribution $F_D$ if \newline $k \in \{a_n : n \geq 1 \}$, and according
to the projection $F_d$ of $F_D$ to the first $d$ coordinates if 
\newline $k \notin \{ a_n \}$.  We assume that:
\begin{eqnarray} \label{eq22}
&& F_D \mbox{ makes the $D$ coordinates independent, and}\\[1ex]
&& F_D \mbox{ has mean zero and finite variance.} \nonumber
\end{eqnarray}

We first state an easy qualitative proposition which is sharpened in
Theorem~\ref{th4.4} below.  
\begin{pr} \label{pr4.2}
Fix distributions $F_2$ and $F_3$ satisfying (\ref{eq22}). If the sequence
$\{ a_n\}$ grow sufficiently fast, then the resulting $\Z^2$ in $\Z^3$
random walk in varying dimension is recurrent.
\end{pr}

\noindent{\sc Proof:}  Denote by $\pi_z$ projection to the $z$-axis
and by $\pi_{xy}$ the projection map to the $x$-$y$ plane.  Since 
$\{ \pi_{xy} (S_k) \}$ is a recurrent planar random walk, we may 
select $a_n$ inductively to satisfy 
\begin{equation} \label{eq23}
P [ \exists k \in (a_n , a_{n+1}] : \pi_{xy} (S_k) = 0] \geq 1/2 .
\end{equation}
The process $\{ \pi_z (S_{a_n}) \}$ is a recurrent one-dimensional
random walk, so there is almost surely a random infinite sequence
$N(1) , N(2), \ldots$ for which $\pi_z (S_{a_{N(j)}}) = 0$ for all 
$j \geq 1$.  Now condition on the sequence $\{\pi_z (S_k) \}$, and use
independence of $\{ \pi_{xy} (S_k) \}$ and
$\{ \pi_z (S_k) \}$ to conclude the following: \newline
With probability at least 1/2,
there are infinitely many $j$ for which
there exists a time \newline
 $k \in (a_{N(j)} , a_{N(j)+1}]$ such that $S_k = 0$.
By the zero-one law, this proves recurrence.   $\Cox$

The argument above is quite general and extends in an obvious way
to the product of two recurrent Markov chains.  Iterating this argument
yields the next corollary.
\begin{sloppypar}
\begin{cor} \label{cor4.3}
If $d_1 < d_2 < \cdots < d_N$ and 
\begin{equation} \label{eq24}
\max \{d_{j+1} - d_j : 1 \leq j \leq N-1 \} \leq 2 \, ,
\end{equation}
then there exists a recurrent process $\{ S_k \}$
with independent increments, which interlaces infinitely many
$d_j$-dimensional steps for each $j$.  More precisely, 
$S_{k+1} - S_k$ has a truly $D(k)$-dimensional distribution
for each $k$, and the sequence $\{ D(k) \}$ takes on only the
values $d_1 , \ldots , d_N$, each one infinitely often.  
If~(\ref{eq24}) is violated then  any such process 
$\{ S_k \}$ must be transient.
\end{cor}
\noindent (To justify the last assertion, observe that the projection
of $\{S_k\}$ to the $d_j+1,...,d_{j+1}$ coordinates
performs a  random walk
in $d_{j+1} - d_j $ dimensional space (up to a time change),
 which is necessarily transient if $d_{j+1}-d_j > 2 $ . ) 
\end{sloppypar}

Next we state the quantitative version, Theorem~\ref{th4.4},
of Proposition~\ref{pr4.2}.  This will be proved in detail.  
We also state similar theorems for RWVD in 2 and 4 dimensions 
and RWVD in 1 and 3 dimensions and give the necessary modifications
to the proof of Theorem~\ref{th4.4}.  Define 
\begin{eqnarray}
\phi (n) & = &  {\log (a_{n+1} / a_n) \over \log a_{n+1}} ; 
   \label{eqdef1} \\[2ex]
\phi_1 (n) & = &  \sqrt {a_{n+1} - a_n \over a_{n+1}} .
   \label{eqdef2} 
\end{eqnarray}

\begin{th} \label{th4.4}
For the $\Z^2$ in $\Z^3$ random walk in varying dimension $\{ S_k \}$
considered in Proposition~\ref{pr4.2}, we have:
\begin{quote}
$(i)$ If $\sum_n n^{-1/2} \phi (n) < \infty$ then $\{ S_k \}$ is
transient.

$(ii)$ If $\sum_n n^{-1/2} \phi (n) = \infty$ and the sequence 
$\{ \phi (n) \}$ is nonincreasing, then $\{ S_k \}$ is recurrent.
\end{quote}
\end{th}

\noindent{\bf Remarks:}  
\begin{enumerate}
\item In particular, $S_k$ is recurrent for
$a_n = \exp (e^{n^{1/2}})$ and transient for $a_n = 
\exp (e^{n^\theta})$ when $\theta  < 1/2$.  

\item The monotonicity assumption in $(ii)$ is far
{}from necessary, and may be weakened in several ways.  If 
$\phi$ is bounded below, $\{S_k\}$ is recurrent and the proof
is easier.  If 
\begin{equation} \label{eqALT}
\sup_{m > n} \phi (m) / \phi (n) < \infty ,
\end{equation} 
then $\{S_k\}$ is still recurrent when $\sum_n n^{-1/2} \phi (n)
= \infty$.  On the other hand, the hypothesis may not be 
discarded completely.  To see this, let $A \subseteq 
\{ 1 , 2 , 3 , \ldots \}$ be a set of times such that
a simple random walk $\{ Y_n \}$ on $\Z^1$ will have
$Y_n = 0$ for only finitely many $n \in A$ almost surely,
even though $\sum_{n \in A} \P [Y_n = 0] = \infty$ 
(cf.\ Example~1 in Section~2 of Benjamini, Pemantle and Peres (1994).)
Define the sequence $\{ a_n \}$ by $a_{n+1} = 2 a_n - 1$ if
$n \notin A$ and $a_{n+1} = a_n^2$ if $n \in A$.  Each
$n \in A$ satisfies $\phi (n) = 1/2$, so the sum in~$(ii)$ is
infinite by the assumption $\sum_{n \in A} \P [Y_n = 0] = \infty$
and the local CLT.
But with probability one,  $S_{a_n}$ is in the $x$-$y$ plane 
for only finitely many $n \in A$, while by Lemma~\ref{lem4.5} , 
$\{ S_k \}$ visits the origin finitely often in time intervals
$[a_n , a_{n+1}]$ for $n \notin A$.  

\end{enumerate}
   
\begin{th} \label{th4.4b}
For any ``$\Z^2$ in $\Z^4$'' random walk in varying dimension $\{S_k\}$,
where  the increment distributions satisfy (\ref{eq22}) : \newline
\noindent{$(i)$} If $\sum_n n^{-1} \phi (n) < \infty$ then $\{ S_k \}$ is
transient.

\noindent{$(ii)$} If $\sum_n n^{-1} \phi (n) = \infty$ and the sequence 
$\{ \phi (n) \}$ is nonincreasing, then $\{ S_k \}$ is recurrent.
\end{th}
(In particular, this inhomogenous walk is recurrent if $a_n=\exp ( e^n ) $, and
 transient if \linebreak
$a_n= \exp  ( e^{n^{\theta }} ) $ \, with $\theta < 1$.) 

\begin{th} \label{th4.4c}
For any ``$\Z^1$ in $\Z^3$'' random walk in varying dimension,
if the increment distributions satisfy (\ref{eq22}) then \newline
\noindent{$(i)$} If $\sum_n n^{-1} \phi_1 (n) < \infty$ then $\{ S_k \}$ is
transient.

\noindent{$(ii)$} If $\sum_n n^{-1} \phi_1 (n) = \infty$ and the sequence 
$\{ \phi_1 (n) \}$ is nonincreasing, then $\{ S_k \}$ is recurrent.
\end{th}
(In particular, this inhomogenous walk  is recurrent if 
$a_n=\exp ( n / {\log ^2 ( n ) } )  $
, and
 transient if the exponent $2$ in the last formula is replaced by
any larger exponent. )

\section{Proofs}

The proofs begin with some elementary estimates on the probability
of returning to the origin in a specified time interval.  
\begin{lem} \label{lem1D}
Let $\{ S_k \}$ be the partial sums of an aperiodic random walk on the 
one-dimensional integer lattice with mean zero and finite 
variance.  Then there exist constants $c_1$ and $c_2$ depending
only on the distribution of the increments, such that
for sufficiently large integers $0 < a < b$,
\begin{equation} \label{eq1D}
c_1 \sqrt{b - a \over b} \leq \P [ S_k = 0 \mbox{ for some } a \leq k < b ]
  \leq c_2 \sqrt{b - a \over b} .
\end{equation}
\end{lem}
\begin{lem} \label{lem2D}
Let $\{ S_k \}$ be the partial sums of an aperiodic random walk on the 
two-dimensional integer lattice with mean zero and finite 
variance.  Then there exist constants $c_1$ and $c_2$ depending
only on the distribution of the increments, such that
for sufficiently large integers $0 < a < b$,
\begin{equation} \label{eq2Da}
c_1 {\log (b/a) \over \log b} \leq \P [ S_k = 0 \mbox{ for some } 
   a \leq k < b ] ,
\end{equation}
and, in the case that $b > 2a$,
\begin{equation} \label{eq2Db}
\P [ S_k = 0 \mbox{ for some } a \leq k < b ] \leq c_2 {\log (b/a) 
   \over \log b} .
\end{equation}
\end{lem}

\noindent{\sc Proof of Lemma~\ref{lem1D}:}  The Local Central 
Limit Theorem (cf. Spitzer (1964)) gives 
\begin{equation} \label{eqLCLT1}
\P [ S_k = 0 ] = {c \over \sqrt{k}} (1 + o(1))
\end{equation}
for some constant $c$ as $k \rightarrow \infty$.  Write
$G$ for the event that $S_k = 0$ for some $k \in [a , b-1]$.  Then
\begin{equation} \label{eqren}
 \P [G] = { \E \# \{ k : a \leq k < b \mbox{ and } S_k = 0 \} \over
    \E \left ( \# \{ k : a \leq k < b \mbox{ and } S_k = 0 \} \| 
    G \right )}  \; . 
\end{equation}
Using (\ref{eqLCLT1}) shows that as $a \rightarrow \infty$,
the  numerator is
$$(c+o(1)) (\sqrt{b} - \sqrt{a}) \geq (c+o(1))(b-a)/{(2 \sqrt{b})}.
$$
 To get an upper bound on the denominator in (\ref{eqren}),
let $T = \min \{ a \leq k < b : S_k = 0 \}$ be the (possibly infinite) 
hitting time and condition on $T$ to get
\begin{eqnarray*}
&& \E \left ( \# \{ k : a \leq k < b \mbox{ and } S_k = 0 \} \| G \right )
 \\[2ex]
& \leq & \sup_{a \leq t < b} \E \left ( \# \{ k : a \leq k < b 
   \mbox{ and } S_k = 0 \} \| T = t \right ) \\[2ex]
& = &  \E \# \{ k : 0 \leq k < b-a 
   \mbox{ and } S_k = 0 \} 
\\[2ex] & \leq &
 C \sqrt{b-a} , 
\end{eqnarray*}
 for some constant $C$ and all positive integers $a<b$.  Thus
$$\P [G] \geq {(c + o(1)) (b-a) / (2 \sqrt{b}) \over C 
   \sqrt{b - a}} \geq c_1 \sqrt{b-a \over b} \, ,$$
for some constant $c_1$ and all sufficiently large $a$.

To prove the second inequality, recompute
$$ \P [G] = { \E \# \{ k : a \leq k < 2b-a \mbox{ and } S_k = 0 \} \over
    \E \left ( \# \{ k : a \leq k < 2b-a \mbox{ and } S_k = 0 \} \| 
    G \right )}  \; . $$
The numerator is now $(c + o(1)) (\sqrt{2b-a} - \sqrt{a})$
which is at most $(2c + o(1)) (b - a) / \sqrt{a}$.
The denominator is at least
$$\inf_{a \leq t < b} 
    \E \left ( \# \{ k : a \leq k < 2b-a \mbox{ and } S_k = 0 \} \| 
    T = t \right ) \geq (c + o(1)) \sqrt{b-a} . $$
Taking the quotient proves the second inequality in the case 
$b \leq 2a$; the case $b > 2a$ is trivial.  
 $\Cox$

\noindent{\sc Proof of Lemma}~\ref{lem2D}:  The Local CLT now gives
$$ \P [ S_k = 0 ] = {c \over k} (1 + o(1)) .$$
Defining $G$ and $T$ as in the preceding proof, it again follows that
$$ \P [G] = { \E \# \{ k : a \leq k < b \mbox{ and } S_k = 0 \} \over
    \E \left ( \# \{ k : a \leq k < b \mbox{ and } S_k = 0 \} \| 
    G \right )}  \; . $$
Using the Local CLT and conditioning on $T$ as before, shows this
to be at least
$$(1 + o(1)) {\log b - \log a \over \log (b-a)} , $$
which proves~(\ref{eq2Da}). 
On the other hand, using the alternate expression
$$ \P [G] = { \E \# \{ k : a \leq k < 2b-a \mbox{ and } S_k = 0 \} \over
    \E \left ( \# \{ k : a \leq k < 2b-a \mbox{ and } S_k = 0 \} \| 
    G \right )} $$
gives
$$\P [G] \leq (1 + o(1)) { \log (2b-a) - \log a \over \log (b-a)}$$
which is at most $c_2 { \log (b/a) \over \log b} $ as long as $b > 2a$,
proving~(\ref{eq2Db}).     $\Cox$

\noindent{\sc Proof of Theorem~\ref{th4.4}:}  The second moment
method will be used.
It is possible to get a good second moment estimate on
the number of intervals $[a_n , a_{n+1} - 1]$ that contain
a return to zero, but only after pruning the short intervals.
We must first prove:
\begin{lem} \label{lem4.5} 
The number of $k$ for which $S_k = 0$ and $a_n \leq k < a_{n+1}$
for some $n$ satisfying $a_{n+1} < 2 a_n$ is almost surely finite.
\end{lem}

\noindent{\sc Proof:}  Let $m(1), m(2) , \ldots$ enumerate
the integers $m$ for which $[2^{m-1} , 2^{m+1} - 1]$ contains 
some $a_n$.  It suffices to show that
finitely many intervals of the form $[2^{m(j)-1} , 2^{m(j)+1} - 1]$ 
contain values of $k$ for
which $S_k = 0$, since these cover all intervals of the
form $[a_n , a_{n+1} - 1]$ satisfying $a_{n+1} < 2 a_n$.  

Fix $j$ and let $n(j)$ denote the least $n$ such that
$a_n \in [2^{m(j) - 1} , 2^{m(j) + 1} - 1]$.  By the independence 
of the coordinates
of $\{ S_k \}$, and by the Local CLT in one and two dimensions,
one sees that for each $k \in [2^{m(j) - 1} , 2^{m(j) + 1} - 1]$, 
the probability of $S_k = 0$ is at most $c / (k \sqrt{n(j)})$.
Summing this over all $k$ in the interval gives
$$\P [S_k = 0 \mbox{ for some } 2^{m(j) - 1} \leq k < 2^{m(j) + 1}]
   \leq {c \over \sqrt{n(j)}} \, . $$
Another way to get an upper bound on this is to see that 
the probability of this event is at most the product of the
probability that the walk returns to the $x$-$y$ plane during
the interval with the probability that it returns to the $z$-axis
during the interval.  Lemmas~\ref{lem1D} and~\ref{lem2D} applied
to the intervals $[n(j) , n(j+1) - 1]$ and $[2^{m(j)-1} , 2^{m(j)+1}
- 1]$ respectively show this product to be at most
$$ c \sqrt{n(j+1) - n(j) \over n(j+1)} \, {1 \over m(j)} .$$
Since $m(j) \geq j$, these two upper bounds may be written as
$$\P [S_k = 0 \mbox{ for some } 2^{m(j) - 1} \leq k < 2^{m(j) + 1}]
   \leq c \min \left ( {1 \over \sqrt{n(j)}} \, , \, \sqrt{n(j+1) - n(j)
   \over n(j+1)} \, {1 \over m(j)} \right ) .$$
Lemma~\ref{lem4.6} with $b_j = n(j+1) - n(j)$ now shows that these
probabilities are summable in $j$, and Borel-Cantelli finishes
the proof.   For continuity's sake, the lemma (which is a fact
about deterministic integer sequences) is given at the
end of the section.    $\Cox$ 

\noindent{\sc Proof of Theorem~\ref{th4.4}} (continued):  Let
$I_n = 1$ if $a_{n+1} \geq 2 a_n$ and $S_k = 0$ for some $k \in
[a_n , a_{n+1} - 1]$, and let $I_n = 0$ otherwise.  Part $(i)$ of
the theorem is just Borel-Cantelli: the hypothesis in $(i)$ and
the estimate~(\ref{eq2Db}) in the case $b > 2a$ together imply
that $\E I_n \leq n^{-1/2} \phi (n)$ is summable.  Thus the
random walk visits zero finitely often in intervals $[a_n , 
a_{n+1} - 1]$ for which $a_{n+1} \geq 2 a_n$; this, together with
Lemma~\ref{lem4.5}, proves $(i)$.

To prove $(ii)$, it suffices, by the 0-1 law (Lemma~\ref{lem4.1}),
to show that the probability of $S_k$ returning to the origin 
infinitely often is positive.  This follows from the
two assertions: $\sum_{n=1}^\infty \E I_n = \infty$, and
$\E (\sum_{n=1}^M I_n)^2 \leq c (\sum_{n=1}^M \E I_n)^2$ 
(cf. Kochen and Stone (1964)).

Seeing that $\sum_{n=1}^\infty \E I_n = \infty$ is easy,
since $\sum_{n=1}^\infty n^{-1/2} \phi (n)$ is assumed
to be infinite; the difference between the two sums is
$$\sum_{n=1}^\infty {\phi (n) \over n^{1/2}} \one_{\{a_{n+1} < 2 a_n\}} .$$
Letting $m(j)$ enumerate those integers $m$ such that $2^{m-1} \leq a_n <
2^{m+1}$ for some $n$, the difference comes out to at most
$$\sum_{j=1}^\infty {\log 4 \over j^{1/2} \log (2^{m(j) - 1})} $$
which is summable since $m(j) \geq j$.  For use below, let
$C_0$ denote this finite sum.

For the second moment computation, take $M$ large enough so that
$\sum_{n=1}^M \E I_n \geq 1$.
  The expected square of $\sum_{n=1}^M  I_n$ 
may be expanded into terms $\E I_n I_r$,  which
we now bound.  Of course, $\E I_n I_r = 0$ if $a_{n+1} < 2 a_n$ 
or $a_{r+1} < 2 a_r$.  Assume now that $\E I_n I_r > 0$ and
that $n < r$ are not consecutive among numbers $k$ with
$\E I_k > 0$.  Then
\begin{eqnarray*}
\E I_n I_r & = & (\E I_n) \E (I_r \| I_n = 1) \\[2ex]
& \leq & c {\phi (n) \over n^{1/2}} \max_{t < a_{n+1}}
   \E (I_r \| S_t = 0) .
\end{eqnarray*}
The inequality~(\ref{eq2Db}) may be used with $b = a_{r+1} - t$
and $a = a_r - t$ to see that
$$\E I_n I_r \leq c \max_{t < a_{n+1}} {\phi (n) \over n^{1/2}}  
   {\log (a_{r+1} - t) - \log (a_r - t) \over \sqrt{r-n} \log (a_{r+1} - t)} 
   .$$
As $t$ increases, the numerator of the last term increases and
the denominator decreases, and since $t < a_{n+1} < a_r/2$, this yields
\begin{eqnarray*}
\E I_n I_r & \leq & c {\phi (n) \over n^{1/2}} {\log (a_{r+1} - a_r/2) 
   - \log (a_r - a_r/2) \over \sqrt{r-n} \log (a_{r+1} - a_r/2)}  \\[2ex]
& \leq & c' {\phi (n) \over n^{1/2}} {\phi (r) \over (r-n)^{1/2}} 
 \leq  c' {\phi (n) \over n^{1/2}} {\phi (r-n) \over (r-n)^{1/2}} \, ,
\end{eqnarray*}
since $\phi (n)$ is assumed to be monotone decreasing.  

Using the bound $\E I_n I_r \leq \E I_n$ for consecutive or identical
nonzero terms, and changing variables $l=r-n$, yields
$$
      \E (\sum_{n=1}^M I_n)^2  \leq  5 c' (C_0 + \sum_{n=1}^M \E I_n)^2.  
$$
By our choice of $M$, 
the second moment is thus bounded by a constant multiple of the
square of the first moment, completing the proof of the theorem.   $\Cox$

\noindent{\sc Proof of Theorem}~\ref{th4.4b} :
This is completely analogous to the previous proof.
It suffices  to change 
$n^{-1/2} \phi (n)$ to $n^{-1} \phi (n)$ everywhere.  The analogue of 
Lemma~\ref{lem4.5} goes through without alteration, since 
$n^{-1} < n^{-1/2}$, and the second moment estimate is easily completed.

\noindent{\sc Proof of Theorem}~\ref{th4.4c} :  
Here we cannot use a version of
Lemma~\ref{lem4.5} - indeed the critical growth rate of $a_n$
is subexponential.

\noindent{$(i)$} Write $I_n$ for the indicator function of the existence
of a $k \in [a_n , a_{n+1} - 1]$ for which $S_k = 0$. Then
$$c_1 {\phi_1 (n) \over n} \leq \E I_n \leq c_2 {\phi_1 (n) \over n}$$
for all large $n$ by Lemma~\ref{lem1D}.  
Invoking the Borel-Cantelli lemma  proves
part $(i)$.

\noindent{$(ii)$} The assumption of this part ensures that
$ \sum_{n=1}^M {\E I_n} \rightarrow \infty $ as $M \rightarrow \infty \, $ .
 To estimate the second moment,
 we consider seperately the contributions of ``long'' and ``short'' intervals:
\begin{equation} \label{eq4.1}
 \E (\sum_{n=1}^M I_n)^2 = 
 \sum_{n=1}^M \left( \E I_n 
  +  2 \sum_{ r=n+1}^M \E I_n I_r \one_{ \{ a_{r+1} \geq 2 a_{n+1} \} } + 
  2 \sum_{r=n+1}^M \E I_n I_r \one_{\{a_{r+1} < 2 a_{n+1} \}} \right) \, .
\end{equation}
 Each summand in the middle term of~(\ref{eq4.1}) may be estimated using 
 Lemma~\ref{lem1D} :
\begin{equation}
 \E I_n I_r  \leq  
 c {\phi_1 (n) \over n} {\sqrt{(a_{r+1} - a_r) / (a_{r+1} - a_{n+1})}
   \over r - n} 
 \leq  2 c  {\phi_1 (n) \over n} 
   {\phi_1 (r) \over r-n} \, ,
\end{equation}
provided that $a_{r+1} \geq 2 a_{n+1}$.
Monotonicity of $\phi_1 $ allows us to bound $\phi_1 (r)$ from above by
 $\phi_1 (r-n)$. This implies that the middle term of~(\ref{eq4.1}) 
is at most 
$$4 c \sum_{n=1}^M \sum_{r=n+1}^M 
   {\phi_1 (n) \over n} 
   {\phi_1 (r-n) \over r-n} \leq
  4 c (\sum_{n=1}^M \E I_n)^2.  
$$
To bound the last term of~(\ref{eq4.1}), fix $n$ and write
\begin{equation} \label{eqshort}
\sum_{r=n+1}^M \E I_n I_r \one_{\{a_{r+1} < 2 a_{n+1}\}}
 \leq c  {\phi_1 (n) \over n}
   \left ( \sum_{r=n+1}^M {c \over r-n} \sqrt{a_{r+1} - a_r
   \over a_{r+1} - a_{n+1}} \one_{\{a_{r+1} < 2 a_{n+1}\}} \right ) .
\end{equation}
Since $\phi_1 (r)^2 = {a_{r+1} - a_r \over a_{r+1}}$ 
is nonincreasing, and since $a_{r+1} < 2 a_n$ for all nonzero
summands in (\ref{eqshort}), it follows that each difference 
$a_{r+1} - a_{r}$
in these summands is at most twice greater than
any preceding difference $a_{r'+1} - a_{r'}$ where $r' < r$ \, . 
  Thus 
$${a_{r+1} - a_r \over a_{r+1} - a_{n+1}} \leq {2 \over r-n}$$
for each $r$ under consideration, and hence 
(denoting $l=r-n$ in the last step) :
$$ \sum_{r=n+1}^M {1 \over r-n} \sqrt{a_{r+1} - a_r
   \over a_{r+1} - a_{n+1}}  \one_{\{a_{r+1} < 2 a_{n+1}\}}
 \leq  c_0+ \sum_{l=1}^\infty 
   {c \over l^{3/2}} < c'.$$

It now follows from (\ref{eqshort}) and Lemma~\ref{lem1D} that
 the last term of~(\ref{eq4.1}) is at most a constant multiple
of the sum $\sum_{n=1}^M \E I_n$. 
Taking $M$ large enough so that this sum is greater
than $1$,  the second moment bound is established.
$\Cox$

Now that the proofs of Theorems~\ref{th4.4},~\ref{th4.4b} 
and~\ref{th4.4c} are complete, it remains to prove the lemma that was used in
the first proof. 

\begin{lem} \label{lem4.6}
Let $b_1 , b_2 , \ldots$ be any sequence of positive integers and
let $B_n = \sum_{k=1}^n b_k$ be the partial sums.  Then
$$\sum_{n=1}^\infty \min \left ( {1 \over \sqrt{B_{n-1}}} \, , \,
   {1 \over n} \sqrt{b_n \over B_n} \right ) \, < \infty .$$
\end{lem}

\noindent{\sc Proof:}  Breaking down the terms according to 
whether $B_{n-1} \geq n^3$ gives
\begin{eqnarray*}
&& \sum_{n=1}^\infty  \min \left ( {1 \over \sqrt{B_{n-1}}} \, , \,
   {1 \over n} \sqrt{b_n \over B_n} \right ) \\[2ex]
& \leq & \sum_{n=1}^\infty \left ( \one_{\{B_{n-1} \geq n^3\}} {1 \over
   \sqrt{B_n}} + \one_{\{B_{n-1} < n^3\}} {1 \over n} \sqrt{ b_n \over B_n}
   \right ) \\[2ex]
& \leq & C + \sum_{n=1}^\infty \one_{\{B_{n-1} < n^3\}} {1 \over n} 
   \sqrt{b_n \over B_n} .
\end{eqnarray*}
To show that the second term is finite, we estimate the sum over
intervals $[M , 2M]$.  Let $\dd_n = \sqrt{b_n / B_n}$,
and, assuming the sum to be nonzero, let $T = \max 
\{ j \leq 2M : B_j \leq (2M)^3 \}$.  Then
\begin{eqnarray*}
&& \sum_{n=M}^{2M} {1 \over n} \sqrt{b_n \over B_n} \one_{\{B_{n-1} < 
   n^3\}}  \\[2ex]
& \leq & {1 \over M} \sum_{n=M}^{2M} \dd_n \one_{\{B_{n-1} < (2M)^3\}} \\[2ex]
& = & {1 \over M} \sum_{n=M}^{T+1} \dd_n .
\end{eqnarray*}
Since $1 + \dd_n^2 = B_{n+1} / B_n$ and each $\dd_n < 1$, we have
$$\prod_{n=M}^T (1 + \dd_n^2) \leq 2 B_T / B_M \leq 16 M^3.  $$
Taking logs gives
$$\sum_{n=M}^T \dd_n^2 \leq \sum_{n=M}^T (\log 2)^{-1} \log (1
   + \dd_n^2) \leq c \log M$$
for $c > 3 / \log 2$ and large $M$.  By Cauchy-Schwarz,
$$\sum_{n=M}^{T+1} \dd_n \leq  \sqrt{M+2} \sqrt{\sum_{n=M}^T \dd_n^2}
   \leq c \sqrt{ M \log M} .$$
Thus
$$ \sum_{n=M}^{2M} \one_{\{B_{n-1} < n^3\}} {1 \over n} \sqrt{b_n \over B_n}
   \leq c \sqrt{\log M \over M} \, \, .$$
This is summable as $M$ varies over powers of 2, which proves the 
lemma.  $\Cox$

\section{A transient inhomogenous random walk with fair bounded steps 
in one and two dimensions}

Let $\{ a_n \}$ and $\{ b_n \}$ be sequences of positive integers,
and let $\{ S_k \}$ be an inhomogenous random walk
 in $\Z^2$ which does $a_1$ steps
of a random walk uniform over all four diagonal neighbors, then $b_1$ 
steps of a horizontal simple random walk, then $a_2$ diagonal steps, 
then $b_2$ horizontal steps, and so on.  

\begin{pr} \label{pr6.1}
If the sequences $\{ a_n \}$ and $\{ b_n \}$ satisfy 
conditions~(\ref{eq dumb}) and~(\ref{eq.t})
below then  $\{ S_k \}$ is transient.
(In particular these conditions hold for $a_n=n^2$ and $b_n=2^n$ .)
\end{pr}

\noindent{\sc Proof:}  Denote 
$A_n = \sum_{j=1}^n a_j$  and    $B_n = \sum_{j=1}^n b_j \, . $
Let $(X_k , Y_k)$ denote the coordinates of $S_k$, and observe that the two
sequences $\{ X_k\}$ and $\{Y_k\}$ are independent of each other.
The first sequence is a simple random walk on the integers, while 
the increments $Y_{k+1}-Y_k$
are simple RW steps if $k \in [A_n+B_n, A_{n+1}+B_n)$ for some $n$;
for $k$ in the complementary time intervals, $Y_{k+1} = Y_k$.

If  $A_n + B_{n-1} \leq k < A_n + B_n$ then $Y_k =Y_{A_n+B_{n-1}}$
 is the position
of a simple random walk after $A_n$ steps, so $Y_k$ vanishes with probability
$1/ {\sqrt{A_n}}$, up to a bounded factor.
The probability that $X_k = 0$ for some $k \in [A_n+B_{n-1}, A_n+B_n)$
may be estimated by \newline
$C [\sqrt{A_n + B_n} - \sqrt{A_n + B_{n-1}}] / \sqrt{b_n} \, $
 (using Lemma \ref{lem1D}). Thus the summability condition
\begin{equation} \label{eq dumb}
\sum_n  {{1 \over \sqrt{A_n}} { \sqrt{A_n + B_n} - \sqrt{A_n + B_{n-1}} \over
   \sqrt{b_n}}} < \infty
\end{equation}
rules out infinitely many returns of $S_k$ to the origin in intervals \newline 
$k \in [A_n + B_{n-1}, A_n + B_n) \,$. 

The only other way $\{ S_k \}$ can return to the origin
 is  during intervals of the type \newline
$[A_n + B_n, A_{n+1} + B_n)$.  The analogous calculation
yields the condition
\begin{equation} \label{eq.t}
\sum_{n=1}^\infty {1 \over \log a_n}
\sum_{j=1}^{a_n} (B_n + A_n + j)^{-1/2} (A_n + j)^{-1/2} 
    \, < \infty
\end{equation}
which is sufficient to ensure finitely many returns 
of $S_k$ to the origin during such intervals.
   $\Cox$

To see a concrete example where Proposition \ref{pr6.1} applies,
with $\{ a_n \}$ and $\{ b_n \}$ growing
almost as slowly as is allowed, let $a_n = (\log n)^{2 + \ee}$
and $b_n = (\log n)^{4 + \ee}$ for some $\ee > 0$.

\noindent{\bf Remark:}  Durrett, Kesten and Lawler (1991)
analyze a random walk in one dimension that interlaces several increment
distributions all having mean zero.  In that setting, distributions
without second moments are necessary in order to obtain transience.


\begin{thebibliography}{YMN}

\bibitem{BPP}
Benjamini, I., Pemantle, R. and Peres, Y. (1994).  Martin capacity for
Markov chains.  {\em Preprint.}

\bibitem{DKL}
Durrett, R., Kesten, H. and Lawler, G. (1991).  Making money in fair
games.  In: Random walks, particle systems and percolation, Durrett
and Kesten Eds.  Birkh\"auser: New York.

\bibitem{KS}
Kochen, S. and Stone, C. (1964). A note on the Borel-Cantelli Lemma.
{\em Illinois J. of Math.} {\bf 8}, 248--251.

\bibitem{Ly}
Lyons, T. (1983). Transience of reversible Markov chains.
{\em Ann.\ Probab.\ } {\bf 11}, 393--402.

\bibitem{Sc}
Scott, D. (1990).  A non-integral-dimensional random walk.  {\em J. Theor.
Prob.} {\bf 3}, 1--7.

\bibitem{Sp}
Spitzer, F. (1964).  Principles of random walk.  Van Nostrand: New York.


\end{thebibliography}
\end{document}